\documentclass[proc]{edpsmath}
\usepackage[utf8]{inputenc}
\usepackage[english]{babel}
\usepackage{amsmath}
\usepackage{mathrsfs}
\usepackage{amsthm}
\usepackage{enumerate}
\usepackage{color}
\usepackage{stmaryrd}
\usepackage{graphicx}
\usepackage[lmargin=1.1in,rmargin=1.1in,tmargin=1.1in,headsep=.1in]{geometry}

\usepackage{verbatim,euscript,eufrak}

\newcommand{\RR}{\mathbb R}

\newcommand{\pat}{\partial_t}
\newcommand{\pax}{\partial_x}

\begin{document}
	
	\title{Modelling and simulation of a wave energy converter}
	%
	\author{Edoardo Bocchi}\address{Departamento de Análisis Matemático \& Instituto de Matemáticas de la Universidad de Sevilla (IMUS), Universidad de Sevilla,
		Avenida Reina Mercedes, 41012 Sevilla, Espa\~{n}a (\email{ebocchi@us.es})}
	\author{Jiao He}\address{Laboratoire de Mathématiques et Modélisation d'Evry (LaMME), Université d'Evry Val d'Essonne,
		23 Boulevard de France, 91037, Evry Cedex, France (\email{jiao.he@univ-evry.fr})}
	\author{Gastón Vergara-Hermosilla}\address{Institut de Math\'ematiques de Bordeaux (IMB), Universit\'e de Bordeaux,
		351 Cours de la Lib\'eration, 33405 Talence Cedex, France (\email{coibungo@gmail.com})}

	\begin{abstract} 
In this work we present the mathematical model and simulations of a particular 
wave energy converter, 
the so-called oscillating water column. In this device, waves governed by the one-dimensional nonlinear shallow water equations arrive from offshore, encounter a step in the bottom and then arrive into a chamber to change the volume of the air to activate the turbine. The system is reformulated as two transmission problems: one is related to the wave motion over the stepped topography and the other one is related to the wave-structure interaction at the entrance of the chamber. We finally use the characteristic equations of Riemann invariants to obtain the discretized transmission conditions and we implement the Lax-Friedrichs scheme to get numerical solutions.

	\end{abstract}
	%
	%
	%
	\maketitle
	\section{Introduction}
	\subsection{General setting}
\begin{figure}
		\centering
		\includegraphics[width=0.9\textwidth]{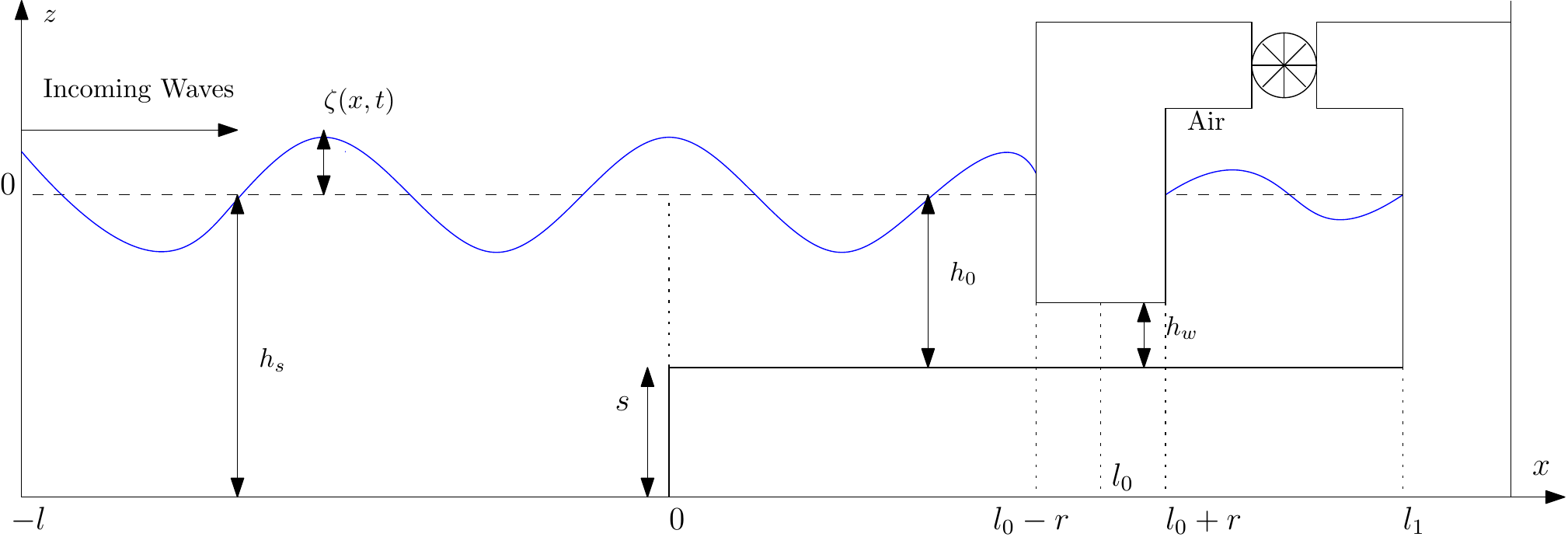}
		\caption{Configuration of the OWC}
		\label{OWC}
	\end{figure}

	This work is devoted to model and simulate an on-shore oscillating water column (OWC), which is a particular type of wave energy converter (WEC) that transforms the energy of waves reaching the shore into electric energy. The structure is installed at the shore in such a way that the water partially fulfills a chamber, which is connected with the outside through a hole where a turbine is placed (see Figure \ref{OWC}). Incoming waves collide with the exterior part of the immersed wall and, after the collision, one part of the wave is reflected while the other part passes below the fixed partially immersed wall and enters the chamber. This increases the water volume inside the chamber and consequently, it creates an airflow that actives the turbine by passing through it and the same occurs when the volume of water reduces inside the chamber. The perpetuation of the incoming waves makes the water inside the chamber oscillate and act as a liquid piston, whose oscillations create electric energy. In this work the wave energy converter is deployed with stepped bottom, which means that
	incoming waves encounter a step in the bottom topography just before reaching the structure. The influence of such step in the OWC device will be discussed later in Section 4. The present research is essentially motivated by a series of works by Rezanejad and collaborators on the experimental and numerical study of nearshore OWCs, in particular, we refer to Rezanejad and Soares \cite{RezSoa18}, where the authors used a linear potential theory to do simulations and showed the improvement of the efficiency when a step is added. Our goal is to numerically study this type of WEC considering as the governing equations for this wave-structure interaction the nonlinear shallow water equations derived by Lannes in \cite{Lan17}, whose local well-posedness was obtained by Iguchi and Lannes in \cite{IguLan18} in the one-dimensional case and by Bocchi in \cite{Boc19} in the two-dimensional axisymmetric case. In the Boussinesq regime and for a fixed partially immersed solid similar equations were studied by Bresch, Lannes and M\'etivier in \cite{BreLanMet19} and in the shallow water viscous case by Maity, San Mart{\'i}n, Takahashi and Tucsnak in \cite{MaiSMTakTuc19} and by Vergara-Hermosilla, Matignon, and Tucsnak in \cite{VH19}.

	We consider an incompressible, irrotational, inviscid and homogeneous fluid in a shallow water regime, which occurs in the region where the OWC is installed. Following \cite{Lan17}, the motion of the fluid is governed by the 1D nonlinear shallow water equations 
	\begin{equation}\label{shalloweq}
	\begin{cases}
	\begin{aligned}
	&\pat \zeta + \pax q=  0\\
	&\pat q + \pax \left(\frac{q^2}{h}\right) + gh \pax \zeta = -\frac{h}{\rho}\pax\underline{P}
	\end{aligned}
	\end{cases}
	\qquad\mbox{for}\quad  x\in (-l, l_1),
	\vspace{0.5em}
	\end{equation}
	where $\zeta(t,x)$ is free surface elevation, $h(t,x)$ is the fluid height, $\rho$ is the fluid density, $\underline{P}$ is the surface pressure of the fluid and $q(t,x)$ is the horizontal discharge defined by $$q(t,x):=\int_{-h_0}^{\zeta(t,x)}u(t,x,z)dz,$$ where $u(t,x,z)$ is the horizontal component of the fluid velocity vector field.\\ 
Let us first give the boundary conditions related to \eqref{shalloweq}. The relevance of these boundary conditions will be explained in Section 2.2. The boundary conditions on the horizontal discharge are  
	\begin{equation}\label{qbc}
	\begin{aligned}
	q \quad \mbox{is continuous} \quad &\mbox{at}\quad  x=0, \, x=l_0\pm r,\\
	q=0 \quad &\mbox{at}\quad  x=l_1,
	\end{aligned}
	\end{equation}and the boundary conditions on the surface elevation are 
	\begin{equation}\label{zetabc}
	\begin{aligned}
	\zeta = f\quad &\mbox{at}\quad  x=-l,\\
	\zeta\quad \mbox{is continuous} \quad &\mbox{at}\quad  x=0,
	\end{aligned}
	\end{equation}
	where $f$ is a prescribed function depending only on time. The surface pressure is given by the constant atmospheric pressure where the fluid is directly in contact with the air,  \textit{i.e.}
	\begin{equation}
		\underline{P}=P_{\mathrm{atm}}\quad \mbox{in}\quad (-l, l_0-r)\cup (l_0+r,l_1)\vspace{0.5em}
	\end{equation} 
	and no surface tension is considered here. On the other hand, under the partially immersed structure, the fluid surface elevation is constrained to be equal to the parametrization of the bottom of the solid $\zeta_w$, \textit{i.e.}
	\begin{equation}
		\zeta=\zeta_w \quad \mbox{in}\quad (l_0-r, l_0+r).
	\end{equation}
	To complete the system, we consider an initial configuration where the fluid is at rest, \begin{equation}
	\zeta (0, x)= \begin{cases}
	0 \quad \mbox{in}\quad (-l, l_0-r)\cup (l_0+r,l_1)\\
	\zeta_w\quad \mbox{in}\quad (l_0-r, l_0+r)
	\end{cases}\quad\mbox{and} \qquad q(0, x)= 0.
	\label{condini}
	\end{equation}

	\subsection{Organization of the paper}
	In Section 2, we derive the model used in the numerical simulations following \cite{IguLan18, Boc19, BreLanMet19}. In particular, we show that the equations \eqref{shalloweq} can be reformulated as two transmission problems, one related to the step in the bottom topography and one related to the wave-structure interaction at the entrance of the chamber.
	Furthermore, the equations in the exterior domain are written as two transport equations on Riemann invariants. In Section 3, we discretize the equations in conservative form using the Lax-Friedrichs scheme and use the Riemann invariants to derive the discretization of the entry condition and boundary conditions. In Section 4, we give several computations showing the numerical solutions of the model and compare the OWC device with and without stepped bottom. At the end of this section, we show the accuracy of the numerical scheme to validate our computations and we discuss the absorbed power and the efficiency of the OWC.

	\subsubsection{Notations}\label{notations}
	We divide the domain of the problem $(-l, l_1)$ into two parts. The interval $\mathcal{I}=(l_0-r,l_0+r)$ is called {\it interior domain}, which is the projection onto the line of the wetted part of the structure, and its complement $\mathcal{E}=(-l, l_1 )\setminus \overline{\mathcal{I}}$, called {\it exterior domain}, which is the union of three intervals $\mathcal{E}_0\cup\mathcal{E}_1\cup \mathcal{E}_2$ with $$\mathcal{E}_0=(-l,0),\quad \mathcal{E}_1=(0, l_0-r) \quad\mbox{ and }\quad \mathcal{E}_2=(l_0+r,l_1),$$
	where $l_1$ is the position of the end of the chamber and  $l_0$ and $r$ are respectively the position of the center and the half length of the partially immersed structure. From the nature of the problem, $l_1>l_0>r$. Moreover, the boundary of $\mathcal{I}$ is formed by the {\it contact points} $\{l_0\pm r\}$, which are the projections on the real line of the triple contact points between fluid, solid and air.  For any function $f$ defined in the real line, its restrictions on the interior domain and the exterior domain are respectively denoted by 
	\begin{equation*}
	f_i:=f_{|_\mathcal{I}}\quad \mbox{and} \quad f_e:=f_{|_\mathcal{E}}.\vspace{0.5em}
	\end{equation*}

	\section{Presentation of the model}\label{derivation}
	\subsection{Governing equations}
	In this section, we present the mathematical model that describes the oscillating water column process considered in this work.
	The model can be essentially divided in three parts: the wave motion over a discontinuous topography represented by the step, the wave-structure interaction at the entrance of the chamber and the wave motion in the chamber. In the exterior domain $\mathcal{E}$, where the fluid is in contact with the air, the surface pressure $\underline{P}_e$ is constrained and is assumed to be equal to the constant atmospheric pressure $P_{\mathrm{atm}}$, while the surface elevation $\zeta_e$ is not known. Contrarily, in the interior domain $\mathcal{I}$, that is the region under the partially immersed structure, the surface elevation $\zeta_i$ is constrained to coincide with the parametrization of the wetted surface, which is assumed to be the graph of some function $\zeta_w$. The surface pressure $\underline{P}_i$ is unknown and it turns out to be a Lagrange multiplier associated with the constraint on $\zeta_i$. For more details on this approach for the study of wave-structure interaction, we refer to \cite{Lan17}. In this work we consider a partially immersed fixed structure with vertical side walls, the parametrization $\zeta_w$ is a constant both in time and space. Summing up, we have an opposite behaviour for the surface elevation and the surface pressure under the structure and elsewhere, that is 
	\begin{align*}
	\zeta_i=\zeta_w, \quad \underline{P}_i \mbox{ is unknown}\qquad &\mbox{and}\qquad  \zeta_e \mbox{ is unknown},\quad\underline{P}_e=P_\mathrm{atm}.
	\end{align*}
	For the exterior domain, we distinguish the region before the step, denoted by $\mathcal{E}_0$ and the region after the step, denoted by $\mathcal{E}_1 \cup \mathcal{E}_2$. The fluid heights are defined respectively by
	\begin{equation*}
	h_e=h_s + \zeta_e \quad \mbox{in}\quad \mathcal{E}_0, \qquad 	h_e=h_0 + \zeta_e \quad \mbox{in}\quad \mathcal{E}_1 \cup\mathcal{E}_2,
	\end{equation*} where $h_s$ and $h_0$ are the fluid heights at rest before the step and after the step respectively. Denoting by $s$ the height of the step, we have $h_s=h_0 +s$.\\
	Therefore the nonlinear shallow water equations \eqref{shalloweq} can be written as the following three systems:
	\begin{enumerate}
		\item for $x\in\mathcal{E}_0,$
		\begin{equation}\label{E0eq}
		\begin{cases}
		\pat \zeta_e + \pax q_e=  0,\\[5pt]
		\pat q_e + \pax \left(\dfrac{q_e^2}{h_e}\right) + gh_e \pax \zeta_e = 0
		\end{cases} \quad \mbox{and} \quad h_e=h_s + \zeta_e,
		\end{equation}
		\item for $x\in\mathcal{E}_1\cup\mathcal{E}_2$
		\begin{equation}\label{E1E2eq}
		\begin{cases}
		\pat \zeta_e + \pax q_e=  0,\\[5pt]
		\pat q_e + \pax \left(\dfrac{q_e^2}{h_e}\right) + gh_e \pax \zeta_e = 0
		\end{cases}\quad \mbox{and} \quad h_e=h_0 + \zeta_e,
		\end{equation}
		
		\item for $x\in\mathcal{I}$,
		\begin{equation}\label{Ieq}
		\begin{cases}
		\pax q_i=  0,\\[5pt]
		\pat q_i = -\dfrac{
			h_w}{\rho}\pax \underline{P}_i
		\end{cases}\quad \mbox{and} \quad h_w=h_0 + \zeta_w.
		\end{equation}
		
	\end{enumerate}

	\subsection{Derivation of the transmission conditions}\label{transm}
	
The  following  section  is  devoted  to  showing  that the motion over the stepped bottom and the wave-structure interaction can  be  reduced  to  two transmission problems for the nonlinear shallow water equations. To do that, we derive the transmission conditions relating the different parts of the model, respectively at the step in the bottom topography and at the side walls of the partially immersed structure.
	
	\subsubsection{At the topography step}
	We consider the problem before the entrance of the chamber not as one shallow water system with a discontinuous topography but rather as a transmission problem between two shallow water systems with flat bottoms where the fluid heights are respectively $h_s +\zeta_e$ and $h_0+\zeta_e$.\\ 
	The first transmission condition is given by the continuity of the surface elevation at the step, namely

	\begin{equation}\label{contiele}
	{\zeta_e}_{|_{x=0^-}}={\zeta_e}_{|_{x=0^+}},
	\end{equation}
	where the traces at $x=0^-$ and at $x=0^+$ are the traces at $x=0$ of the unknowns before the step and after the step respectively.\\ 
	The second transmission condition is given by the continuity of the horizontal discharge at the step, namely
	\begin{equation}
	{q_e}_{|_{x=0^-}}={q_e}_{|_{x=0^+}}.
	\label{continuityq}
	\end{equation}

	\subsubsection{At  the structure side-walls}

	The transmission conditions at the side-walls of the partially immersed structure are derived from the continuity of the horizontal discharge at the side-walls and the assumption that the total fluid-structure energy is equal to the integral in time of the energy flux at the entry of the domain. 
	The continuity of the horizontal discharge at $x=l_0\pm r$ together with the fact that $\pax q_i=0$ gives the first transmission condition between the $\mathcal{E}_1$ and $\mathcal{E}_2$, which reads
	\begin{equation}\label{transqwalls}
	\llbracket q_e \rrbracket:={q_e}_{|_{x=l_0+ r}} - 	{q_e}_{|_{x=l_0- r}}=0.
	\end{equation}
	Let us now derive the second transmission condition at $x=l_0\pm r$. To do that, we show the local conservation of the fluid energy in the exterior domain and in the interior domain as in \cite{BreLanMet19}. \\\\
	{\it Exterior domain.}  Considering the nonlinear shallow water equations in $\mathcal{E}$, multiplying the first equation in \eqref{E0eq}-\eqref{E1E2eq} by 
	$\rho g\zeta_e$  and the second equation by $\rho\dfrac{q_e}{h_e}$, and considering the fact that $\pat h_e=-\pax q_e$, we obtain  
	\begin{equation}\label{eq1enex}
	\left\{
	\begin{aligned}
	&\pat \left( \rho\frac{q_e^2}{2h_e} \right)  +\rho g\zeta_e \pax q_e=  0,\\
	&\pat \left( \rho \frac{q_e^2}{2h_e} \right) - 
	\rho\frac{q_e^2}{2h_e^2}\pax q_e  +\rho\frac{q_e}{h_e} \pax \left( \frac{q_e^2}{h_e}\right)  + \rho gq_e \pax \zeta_e = 0.
	\end{aligned}
	\right.
	\end{equation}
	Adding both equations in (\ref{eq1enex}), we obtain 
	\begin{equation*}
	\pat \left( \rho g
	\frac{\zeta_e^2 }{2}+ \rho\frac{q_e^2}{2h_e}
	\right) 
	+\rho g\zeta_e\pax  q_e+\rho gq_e \pax 
	\zeta_e
	-
	\rho\frac{q_e^2}{2h_e^2}\pax q_e+ \rho\frac{q_e}{h_e} \pax \left( 
	\frac{q_e^2}{h_e}
	\right)
	= 0.
	\end{equation*}
	We compute that 
	\begin{equation*}\label{eq3enex}
	g\zeta\pax  q+gq \pax 
	\zeta = \pax(g\zeta q) \quad \mbox{ and }\quad 
	-
	\frac{q^2}{2h^2}\pax q+ \frac{q}{h} \pax \left(  
	\frac{q^2}{h} \right) = \pax \left( \frac{q^3}{2h^2}\right),
	\end{equation*}
	and, denoting by $\mathfrak{e}_{ext}$ and by $\mathfrak{f}_{ext}$ respectively the local fluid energy and the local flux 
	\begin{equation*}
	\mathfrak{e}_{ext}= \rho\frac{q_e^2}{2h_e} + g\rho \frac{\zeta_e^2}{2} \quad \mbox{ and }\quad  \mathfrak{f}_{ext}=\rho\frac{q_e^3}{2h_e^2} + g\rho\zeta_e q_e,
	\end{equation*}
	we obtain the local conservation of the fluid energy in the exterior domain
	\begin{equation}\label{extcon}
	\pat \mathfrak{e}_{ext} + \pax \mathfrak{f}_{ext}=0.
	\end{equation} 
	{\it Interior domain}. 
	Let us remark that from the first equation in \eqref{Ieq} one gets that  $q_i\equiv q_i(t)$ in the interior domain. Multiplying the second equation in \eqref{Ieq} by $\dfrac{q_i}{h_w}$, we obtain	
	$$
	\pat \left( \rho\frac{ q_i^2}{2h_w}  \right) + \pax \left(q_i \underline{P}_i \right)= 0, 
	$$
	and, denoting by $\mathfrak{e}_{int}$ and by $\mathfrak{f}_{int}$ respectively the local fluid energy and the local flux                       
	\begin{equation*}           
	\mathfrak{e}_{int}=  \rho \frac{q_i^2}{2h_w} + \rho g \frac{\zeta_w^2}{2} \quad \mbox{ and }\quad  \mathfrak{f}_{int}=q_i P_i ,
	\end{equation*}
	we obtain the local conservation of the fluid energy in the interior domain, 
	\begin{equation}\label{intcon}
	\pat \mathfrak{e}_{int} + \pax \mathfrak{f}_{int}=0.
	\end{equation} 
	
	Now we assume that the total fluid-structure energy at time $t$ is equal to the integral between $0$ and $t$ of the sum between the energy flux at the entry of the domain and the difference of the energy fluxes at the step, \textit{i.e.}
	\begin{equation*}
		{E}_{\mathrm{fluid}} + E_{\mathrm{solid}} = \int_{0}^{t}\left({\mathfrak{f}_{ext}}_{|_{x=-l}}+{\mathfrak{f}_{ext}}_{|_{x=0^+}}-{\mathfrak{f}_{ext}}_{|_{x=0^-}}\right),
	\end{equation*} with the fluid energy defined by \begin{equation*}
	{E}_{\mathrm{fluid}} = \int_{\mathcal{I}} \mathfrak{e}_{int} + 
	\int_{\mathcal{E}} \mathfrak{e}_{ext}.
	\end{equation*}This assumption is an adaptation to a bounded domain case of the conservation of total fluid-structure energy assumed in \cite{BreLanMet19}. We remark that the difference of the energy fluxes at the step ${\mathfrak{f}_{ext}}_{|_{x=0^+}}-{\mathfrak{f}_{ext}}_{|_{x=0^-}}$ does not vanish due to the discontinuity of the fluid height at $x=0$ in the presence of the step. The fact that the structure is fixed ($\frac{d}{dt} E_{\mathrm{solid}}=0$) yields
	\begin{equation*}
	\frac{d}{dt}{E}_{\mathrm{fluid}} = \int_{\mathcal{I}} \pat\mathfrak{e}_{int} + 
	\int_{\mathcal{E}} \pat\mathfrak{e}_{ext}={\mathfrak{f}_{ext}}_{|_{x=-l}}+
{\mathfrak{f}_{ext}}_{|_{x=0^+}}-{\mathfrak{f}_{ext}}_{|_{x=0^-}}
	.
	\end{equation*}			
	From (\ref{extcon}) and (\ref{intcon}) we have 
	\begin{equation*}
	-\int_{\mathcal{I}} \pax\mathfrak{f}_{int} 
	-\int_{\mathcal{E}} \pax\mathfrak{f}_{ext}={\mathfrak{f}_{ext}}_{|_{x=-l}}+
	{\mathfrak{f}_{ext}}_{|_{x=0^+}}-{\mathfrak{f}_{ext}}_{|_{x=0^-}}
	.
	\end{equation*}	
	Using the boundary conditions \eqref{qbc} and \eqref{zetabc} we get
	\begin{equation*}
	\llbracket \mathfrak{f}_{int} \rrbracket 
	=
	\llbracket \mathfrak{f}_{ext} \rrbracket, 
	\end{equation*}where the brackets $\llbracket \cdot\rrbracket$ are defined as in \eqref{transqwalls}. By definition of the fluxes it follows
	$$
	\llbracket q_i \underline{P}_i   \rrbracket = \rho \left \llbracket
	\frac{q_e^3}{2h_e^2} + g \zeta_e q_e
	\right \rrbracket 
	$$
	and from \eqref{qbc} and \eqref{transqwalls} we obtain
	\begin{equation*}
	\llbracket \underline{P}_i   \rrbracket = \rho  \left \llbracket
	\frac{q^2_e}{2h_e^2} + g \zeta_e
	\right \rrbracket .
	\end{equation*}
	Integrating on $(l_0+r,l_0-r )$, the second equation in \eqref{Ieq} yields
	\begin{equation*}
	-\frac{\rho \,2r}{h_w} \frac{d}{dt}q_i =	\llbracket \underline{P}_i.   \rrbracket 
	\end{equation*}Combining the last two equalities,
	we get the following transmission condition
	\begin{equation}
	-\frac{\,2r}{h_w} \frac{d}{dt}q_i = \left \llbracket
	\frac{q^2_e}{2h_e^2} + g \zeta_e \right \rrbracket .
	\end{equation}

	\subsection{Reformulation as two transmission problems}	
	
	Coupling the governing equations \eqref{E0eq}-\eqref{Ieq} with the conditions derived in the previous section, we have therefore reduced the problem of the OWC essentially to two transmission problems. The first one in $\mathcal{E}_0 \cup \mathcal{E}_1$  reads:
	\begin{equation}
	\left\{
	\begin{aligned}
	&\pat \zeta_e + \pax q_e=  0,\\
	&\pat q_e + \pax \left(\frac{q_e^2}{h_e}\right) + gh_e \pax \zeta_e = 0,
	\end{aligned}
	\right. \quad h_e=h_s+\zeta_e \quad \mbox{in} \quad \mathcal{E}_0, \quad h_e=h_0+\zeta_e \quad \mbox{in} \quad \mathcal{E}_1, 
	\end{equation}
	with transmission conditions at $x=0$
	\begin{equation}
	{\zeta_e}_{|_{x=0^-}}={\zeta_e}_{|_{x=0^+}}, \qquad {q_e}_{|_{x=0^-}}={q_e}_{|_{x=0^+}}.
	\end{equation}
	The second transmission problem in $\mathcal{E}_1\cup \mathcal{E}_2$ reads:
	\begin{equation}
	\left\{
	\begin{aligned}
	&\pat \zeta_e + \pax q_e=  0,\\
	&\pat q_e + \pax \left(\frac{q_e^2}{h_e}\right) + gh_e \pax \zeta_e = 0,
	\end{aligned}
	\right. \quad h_e=h_0+\zeta_e,
	\end{equation}with transmission conditions at $x=l_0\pm r$
	\begin{equation}\label{transwalls}
	\left \llbracket q \right \rrbracket = 0, \qquad  -\alpha \frac{d}{dt}q_i=\left \llbracket
	\frac{q^2_e}{2h_e^2} + g \zeta_e
	\right \rrbracket,
	\end{equation}
	where $\alpha =\dfrac{2r}{h_w}$ and $h_w=h_0+\zeta_w$.
	
	\subsection{Riemann invariants}\label{Riemann_invariants_section}
	Let us now rewrite the nonlinear shallow water equations \eqref{E0eq} and \eqref{E1E2eq} in the exterior domain $\mathcal{E}$ in a compact form by introducing the couple $U= (\zeta_e, q_e )^T$:
	\begin{equation}\label{compactform}
	\pat U + A(U) \pax U = 0,
	\end{equation}
	where 
	\begin{equation*}
	A(U) = 
	\left(
	\begin{matrix}
	0 & 1 \\
	gh_e - \frac{q_e^2}{h_e^2} & \frac{2q_e}{h_e}
	\end{matrix}
	\right).\vspace{1em}
	\end{equation*}
	The eigenvalues $\lambda_+(U)$ and $-\lambda_-(U)$ of the matrix $A(U)$ and the associated eigenvectors $e_+ (U)$ and $e_- (U)$ are given by 
	$$
	\lambda_+ (U) = \frac{q_e}{h_e} + \sqrt{g h_e }, \quad -\lambda_- (U) = \frac{q_e}{h_e} - \sqrt{g h_e },
	$$
	$$e_+ (U) = \left( \sqrt{gh_e}- \frac{q_e}{h_e}, 1\right)^T ,\quad e_- (U) = \left(-\sqrt{gh_e}- \frac{q_e}{h_e}, 1\right)^T.\vspace{1em}$$
	Notice that $\lambda_+ > 0$ and $\lambda_- > 0$. 
	Taking the scalar product of \eqref{compactform} and eigenvectors, we obtain
	$$
	\pat \left(2 \sqrt{gh_e} \pm \frac{q_e}{h_e}\right) \pm \left(\sqrt{gh_e} \pm \frac{q_e}{h_e}\right) \pax \left(2\sqrt{gh_e} \pm \frac{q_e}{h_e}\right)=0.
	$$
	Let us introduce the right and the left Riemann invariant $R$ and $L$ associated to the nonlinear shallow water equations, respectively
	\begin{equation}\label{RiemannRL}
	R(U):= 2 \left(\sqrt{gh_e} - \sqrt{gh_0}\right) + \frac{q_e}{h_e},\qquad L(U):= 2 \left(\sqrt{gh_e} - \sqrt{gh_0}\right) - \frac{q_e}{h_e}.
	\end{equation}
	Hence we can write the 1D nonlinear shallow water equations in the exterior domain as the two following transport equations on $R$ and $L$: 
	\begin{equation}\label{Riemanneq}
	\pat R(U) + \lambda_+(U) \pax R(U)  = 0, \qquad \pat L(U) - \lambda_-(U) \pax L(U) = 0.
	\end{equation}
	
	We will see that these two transport equations of Riemann invariants are helpful when we solve our model by numerical method. More details about Riemann invariants of the nonlinear shallow water equations can be found in \cite{lannes2019generating}.

\section{Discretization of the model}\label{discret}
		We have reformulated in the previous section the mathematical model of the oscillating water column as two transmission problems. 
		This section is devoted to discretize the nonlinear shallow water equations \eqref{E0eq}-\eqref{Ieq} at the level of the numerical scheme. More precisely, we will use the Lax-Friedrichs scheme to solve our main equations and use Riemann invariants to address the entry conditions and all boundary conditions. 
		
		\subsubsection{Numerical notations}\label{notationsnum}
		We use the following notations throughout this section:
		\begin{itemize}
			\item in our system, the whole numerical domain $[-l,l_0]$ is composed of four parts: $[-l,0]$, $[0,l_0-r]$, $[l_0 -r, l_0 +r]$ and $[l_0+r,l_1]$. Each interval is divided into cells $(\mathcal{A}_i)_{1\leq i \leq n_x}$ with $\mathcal{A}_i = [x_{i-1}, x_i]_{1 \leq i \leq n_x } $ of size $\delta_x$. More precisely, we have 
		\begin{align*}
		    &	x_{0} =-l,  ..., \, x_{i} = -l+i \delta_x ,...,\, x_{n_{1,x}} = 0 ;\\
		    &	x_{n_{1,x}+1} = \delta_x,\, ..., \,x_{n_{1,x}+i} = i \delta_x ,...,\, x_{n_{1,x}+n_{2,x}} = l_0-r ;\\
		   &	x_{n_{1,x}+n_{2,x} +1}=l_0 - r + \delta_x, ..., x_{n_{1,x}+n_{2,x}+i} = l_0 -r + i \delta_x ,..., x_{n_{1,x}+n_{2,x} + n_{3,x}} = l_0+r  ;\\
		   	& x_{n_{1,x}+ n_{2,x} + n_{3,x} +1} =l_0 + r + \delta_x,\, ..., \,x_{n_{1,x}+ n_{2,x} + n_{3,x} +i} = l_0 + r + i \delta_x,\,..., \,x_{n_{1,x}+n_{2,x} + n_{3,x} + n_{4,x}} =l_1, 
		\end{align*}
		with $l=n_{1,x}\delta_x, \, l_0 - r = n_{2,x} \delta_x, \, 2r =  n_{3,x}  \delta_x $ and $l_1 - (l_0+r) = n_{4,x} \delta_x$;
			\item we denote by $\delta_t$ the time step. According to CFL condition, time step $\delta_t$ can be specified by $\delta_x$;
			\item for any quantity $U$, we denote by $U^m_i$ its value at the position $x_i$ at time $t^m = m\delta_t$. For instance, the variables $\zeta_i^m$ denotes the value of the free surface elevation $\zeta$ at the position $x_i$  at time $t^m = m\delta_t$.
		\end{itemize}

	\subsection{Discretization of the equation}\label{discret1}
		The finite difference method is a standard discretization approach for partial differential equations, especially those that arise from conservation laws. We first rewrite equation \eqref{compactform} as the following conservative form : 
		\begin{equation}\label{conver}
		\pat U + \pax( F(U)) = 0 ,
		\end{equation}
		with 
		$$
		F(U) = \left(q_e, \frac{1}{2} g \left(h_e^2 - h_0^2\right) + \frac{q_e^2}{h_e}\right)^T.
		$$
		By means of a finite difference approach, equation \eqref{conver} can be discretized as  
		\begin{equation*}
		\frac{U_{i}^{m+1}-U_{i}^{m}}{\delta_{t}}+\frac{F_{i+1 / 2}^{m}-F_{i-1 / 2}^{m}}{\delta_{x}}=0, 
		\end{equation*}
		where the flux $F$ is discretized with cell centres indexed as $i$ and cell edge fluxes indexed as $i\pm 1/2$. 
		The choice of $F^m_{i\pm 1/2}$ depends on the numerical scheme. We consider here the well-known Lax–Friedrichs scheme proposed by Lax \cite{lax1954weak} to get the discrete flux 
		\begin{equation}
		F_{i-1 / 2}^{m}=\frac{1}{2}\left(F_{i}^{m}+F_{i-1}^{m}\right)-\frac{\delta_{x}}{2 \delta_{t}}\left(U_{i}^{m}-U_{i-1}^{m}\right),
		\end{equation}
		where $i\geq 1$ and $F_i^m=F(U^m_i)$.

		\subsection{Discretization of the entry condition}
		At the entrance of our system, the surface elevation is given by a prescribed function $f$ depending only on time,
		$$\zeta^{m}|_{{x=-l}}=f^{m} :=f\left(t^{m}\right).$$
		In order to express the entry condition for the horizontal discharge, let us first recall that from \eqref{RiemannRL} one has 
		$$
		q_e= \frac{h_e}{2} \left(R-L\right),  \qquad R+L = 4 \left(\sqrt{gh_e} - \sqrt{gh_0}\right),
		$$
		where $R$ and $L$ are respectively the right and the left Riemann invariant associated to the nonlinear shallow water equations. We get
		$$
		q_e = h_e \left(2 \left(\sqrt{gh_e} - \sqrt{gh_0}\right) -L \right).
		$$
		Hence, the value of $q_e$ at $x=-l$ is given by 
		$$
		{q_e}|_{{x=-l}} = \left(h_0+ f(t)\right) \left(2 \left(\sqrt{g (h_0+ f(t))} - \sqrt{gh_0}\right) - L|_{{x=-l}} \right).
		$$
		On the right-hand side of the relation above, $ L|_{{x=-l}}$ is unknown. First we have to determine $ L|_{{x=-l}}$ in order to determine ${q_e}|_{{x=-l}}$. This can be achieved by the transport equation for $L$ in \eqref{Riemanneq}. After discretizing it as in \cite{mar05}, we get
		\begin{equation}
		\frac{L_{0}^{m}-L_{0}^{m-1}}{\delta_{t}}-\lambda_{-} \frac{L_{1}^{m-1}-L_{0}^{m-1}}{\delta_{x}}=0,
		\end{equation}
		where $L_0^m$ is the value of $ L$ at $x=-l$ at time $t^m$ and $\lambda_{-}$ is computed as a linear interpolation between $\lambda_{-,0}$ and $\lambda_{-,1}$ following \cite{lannes2019generating}, namely
		\begin{equation*}
		\lambda_{-} = \beta \lambda_{-,0} + (1-\beta) \lambda_{-,1}
		\end{equation*}
		with $0\leq \beta \leq 1$ such that $\lambda_{-} \delta_t =\beta \delta_x$. Moreover, we can compute $\lambda_{-}$ as
		\begin{equation*}
		\lambda_{-} = \frac{\lambda_{-,1}}{1+ \frac{\delta_t}{\delta_x} \lambda_{-,1} - \frac{\delta_t}{\delta_x} \lambda_{-,0}}.
		\end{equation*}
		Thus, we have
		\begin{equation}
		L_{0}^m = \left(1-\lambda_- \frac{\delta_t}{\delta_x}\right)L_{0}^{m-1}+ \lambda_- \frac{\delta_t}{\delta_x} L_{1}^{m-1},
		\end{equation}
		which gives $L_{0}^m$ in terms of its values at the previous time step and in terms of interior points.   
		
		\subsection{Discretization of the boundary conditions}
		Since our system is composed by four parts, it remains three boundary conditions should be taken into consideration besides the entry condition at $x=-l$. When wave arrives from the offshore, it will encounter a step in the bottom and then arrive into a chamber, and finally arrive to the wall (see the configuration \ref{OWC}). More precisely, the first boundary condition is at the discontinuity of the topography located at $x= 0$ and the second is at the partially immersed structure side-walls located at $x=l_0 \pm r$. The last boundary condition is at the end of the chamber, located at $x=l_1$. 
		\subsubsection{At the topography step}\label{discon}
		Let us first consider the shallow water wave equations with discontinuous topography, namely, it is a system with depth $h_s$ on  $\RR_{-} = \{x<0\}$ and depth $h_0$ on $\RR_+ =  \{x>0\}$. Our equation turns out to be 
		\begin{equation*}
		\pat U + \pax\left( F(U)\right) = 0 ,
		\end{equation*}
		with 
		\begin{equation*}
		F(U) = \left\{
		\begin{aligned}
		&\left(q_e, \frac{1}{2} g \left((h_s + \zeta_e)^2 - h_s^2\right) + \frac{q_e^2}{h_s + \zeta_e}\right)^T, \quad \text{in} \quad (0,T) \times \RR_{-},\\
		&\left(q_e, \frac{1}{2} g \left((h_0 + \zeta_e)^2 - h_0^2\right) + \frac{q_e^2}{h_0 + \zeta_e}\right)^T, \quad \text{in} \quad (0,T) \times \RR_+.
		\end{aligned}
		\right.
		\end{equation*}
		
	From transmission conditions \eqref{contiele} and \eqref{continuityq}, we have the continuity of the surface elevation $\zeta_e$ and of the horizontal discharge $q_e$ at $x=0$:
		\begin{equation}\label{BDCdiscon}
		\zeta_e^l |_{x=0} = \zeta_e^r |_{x=0}, \qquad q_e^l |_{x=0}  = q_e^r |_{x=0}.
		\end{equation}

	Let us denote the right Riemann invariant in the domain $\RR_-$ by $R^l$ and the left Riemann invariant in the domain $\RR_+$ by $L^r$. We then find two expressions of $q_e$ describing $q_e^l |_{x=0}$ and $q_e^r |_{x=0}$, respectively,
		\begin{equation}\label{qdiscon}
		\left\{
		\begin{aligned}
		&q_e^l |_{x=0}  = \left(h_s+ \zeta_e^l |_{x=0}  \right) \left(R^l |_{x=0}  -2 \left(\sqrt{g (h_s+ \zeta_e^l |_{x=0}  )} - \sqrt{g h_s}\right)\right),\\
		& q_e^r  |_{x=0} = \left(h_0 + \zeta_e^r |_{x=0} \right) \left(2 \left(\sqrt{g(h_0+ \zeta_e^r |_{x=0} )} - \sqrt{g h_0}\right) -L^r |_{x=0} \right).
		\end{aligned}
		\right.
		\end{equation}
		According to the relations \eqref{BDCdiscon}, we observe that \eqref{qdiscon} is a system of two nonlinear equations on the two unknowns $\zeta_e^l |_{x=0}$ (respectively $\zeta_e^r |_{x=0}$) and $q_e^l |_{x=0}$ (respectively $q_e^r |_{x=0}$). We write it in the compact form 
		\begin{equation}\label{NLsystem}
		    F(x_1,x_2)=0, 
		\end{equation}
		where $x_1=\zeta_e^l |_{x=0}$, $x_2=q_e^l |_{x=0}$ and the vector $F=(F_1, F_2)$ is given by $$F_1=(h_s+ x_1  ) \left(R^l |_{x=0}  -2 \left(\sqrt{g (h_s+ x_1  )} - \sqrt{g h_s}\right)\right) - x_2,$$
		$$F_2=(h_0 + x_1 ) \left(2 \left(\sqrt{g(h_0+ x_1 )} - \sqrt{g h_0}\right) -L^r |_{x=0} \right)-x_2.\vspace{0.5em}$$
		In the case $h_s=h_0$ (without step) we can derive from \eqref{qdiscon} a third degree equation on $\sqrt{h_0+\zeta_e^l |_{x=0}}$ and take the unique solution that gives $\zeta_e^l |_{x=0}=0$ when $R^l |_{x=0},L^r |_{x=0}=0$ (we refer to \cite{Lan17} for this case).
		Here, since $h_s \neq h_0$, we use MATLAB nonlinear system solver \textit{fsolve} with initial point $(0,0)$ to solve \eqref{NLsystem}. Before doing that, we have to determine the values of the two Riemann invariants $R^l |_{x=0}$ and $L^r |_{x=0}$.
		The transport equations for $R^l$ and $L^r$ are the following: 
		\begin{equation}\label{characteristicdiscon}
		\pat R^l + \lambda_+^l(U) \pax R^l  = 0,\qquad \pat L^r - \lambda_-^r(U) \pax L^r=0,
		\end{equation}
		where the corresponding eigenvalue $\lambda_+^l$ in the domain $\RR_{-}$ is given by 
		\begin{equation}
		\label{lambda_+}
		\lambda_+^l (U) = \frac{q_e}{h_s + \zeta_e} + \sqrt{g (h_s + \zeta_e)},
		\end{equation}
		and the corresponding eigenvalue $-\lambda_-^r$ in the domain $\RR_{+}$ is given by 
		\begin{equation}
		\label{lambda_-}
		-\lambda_-^r (U) = \frac{q_e}{h_0 + \zeta_e} - \sqrt{g (h_0 + \zeta_e) }.
		\end{equation}
		Let us emphasize that we use here the same interpolation for $\lambda_+$ and $\lambda_-$ as in  \cite{mar05}.
		After discretization of equations \eqref{characteristicdiscon}, we get 
		\begin{equation*}
		\frac{(R^l )_{n_{1,x}}^{m}-(R^l)_{n_{1,x}}^{m-1}}{\delta_{t}}+\lambda_{+}^l \frac{(R^l )_{n_{1,x}}^{m-1}-(R^l)_{n_{1,x}-1}^{m-1}}{\delta_{x}}=0, \qquad
		\frac{(L^r)_{n_{1,x}}^{m}-(L^r)_{n_{1,x}}^{m-1}}{\delta_{t}}-\lambda_{-}^r \frac{(L^r)_{n_{1,x}+1}^{m-1}-(L^r)_{n_{1,x}}^{m-1}}{\delta_{x}}=0,
		\end{equation*}
		where $\lambda_{+}^l$, $\lambda_{-}^r $ are as in \eqref{lambda_+}-\eqref{lambda_-} and we recall that $(R^l )_{n_{1,x}}^{m}$ is the value of $R^l$ at $x_{n_{1,x}}$ and $t^m$ (see Notations \ref{notationsnum}).
		Hence, we have
		\begin{equation}\label{discrediscon}
		(R^l )_{n_{1,x}}^{m}=\left(1-\lambda_{+}^l \frac{\delta_{t}}{\delta_{x}}\right) (R^l )_{ n_{1,x}}^{m-1}+\lambda_{+}^l \frac{\delta_{t}}{\delta_{x}} (R^l )_{n_{1,x}-1}^{m-1},\qquad (L^r)_{n_{1,x}}^{m}=\left(1-\lambda_{-}^r \frac{\delta_{t}}{\delta_{x}}\right) (L^r)_{ n_{1,x}}^{m-1}+\lambda_{-}^r \frac{\delta_{t}}{\delta_{x}} (L^r)_{n_{1,x}+1}^{m-1},
		\end{equation}
		which give $(R^l )_{n_{1,x}}^{m}$ and $(L^r)_{n_{1,x}}^{m}$ in terms of their values at the previous time step and in terms of interior points. 
		
		Gathering the relations \eqref{BDCdiscon}, \eqref{qdiscon} and \eqref{discrediscon}, we can solve $\zeta_e^l |_{x=0}$ (respectively $\zeta_e^r |_{x=0}$) and $ q_e^l |_{x=0}$ (respectively $q_e^r |_{x=0}$), which give us the boundary conditions at the step. 
		
		\subsubsection{At  the structure side-walls}
		Compared with the derivation of the boundary conditions near the step, the idea to derive the boundary condition near the fixed  partially immersed structure is almost the same. There are two differences between them. The first one is that, since the depth is always $h_0$, the eq. \eqref{qdiscon} becomes 
		\begin{equation}\label{qobject}
		\begin{aligned}
		& q^l_e |_{x=l_0-r}  = (h_0+ \zeta^l_e |_{x=l_0-r}  ) \left(R^l |_{x=l_0-r}  -2 \left(\sqrt{g (h_0+ \zeta^l_e |_{x=l_0-r}  )} - \sqrt{g h_0}\right)\right),\\
		& q^r_e  |_{x=l_0+r} = (h_0 + \zeta^r_e |_{x=l_0+r} ) \left(2 \left(\sqrt{g(h_0+ \zeta^r_e |_{x=l_0+r} )} - \sqrt{g h_0}\right) -L^r |_{x=l_0+r} \right),
		\end{aligned}
		\end{equation}
		where we denote the horizontal discharge in the exterior domain on the left-hand side of the object by $q^l_e$ and on the right-hand side of the object by $q^r_e$.
		Let us recall that $q_i$ is the horizontal discharge in the interior domain $\mathcal{I}$. From the first transmission condition in \eqref{transwalls}, we know that 
		$$q^l_e |_{x=l_0-r} = q_i =q^r_e |_{x=l_0+r}.$$
		The second difference is that, unlike in the previous subsection, we do not have the continuity condition of $\zeta_e$ at the structure side-walls. Nevertheless, we consider the discretization of the second transmission condition in \eqref{transwalls}, hence we get 
		\begin{equation*}
		-\alpha \frac{(q_e)^m_{l_0-r} - (q_e)^{m-1}_{l_0-r}  }{\delta t}=\frac{\left((q_e^l)^{m-1}_{l_0+r}\right)^2}{2\left(h_0 + (\zeta_e^l)^{m-1}_{l_0+r}\right)^2} +  g (\zeta_e^l)^{m-1}_{l_0+r} -\frac{\left((q_e^r)^{m-1}_{l_0-r}\right)^2}{2\left(h_0 + (\zeta_e^r)^{m-1}_{l_0-r}\right)^2}  - g (\zeta_e^r)^{m-1}_{l_0-r}    
		.
		\end{equation*}where for the sake of clarity $(q_e)^m_{l_0-r}=(q_e)^m_{n_{1,x}+ n_{2,x}}$ and $(q_e)^m_{l_0+r}=(q_e)^m_{n_{1,x}+ n_{2,x}+ n_{3,x}}$ (analogously for $(\zeta_e)^m_{l_0-r}$ and $(\zeta_e)^m_{l_0+r}$).
		Then, $q_e$ at $x=l_0-r$ is expressed as
	\begin{equation}
		(q_e)^m_{l_0-r} = (q_e)^{m-1}_{l_0-r} - \frac{\delta t}{ \alpha} \left( \frac{\left((q_e^l)^{m-1}_{l_0+r}\right)^2}{2\left(h_0 + (\zeta_e^l)^{m-1}_{l_0+r}\right)^2}  - \frac{\left((q_e^r)^{m-1}_{l_0-r}\right)^2}{2\left(h_0 + (\zeta_e^r)^{m-1}_{l_0-r}\right)^2}   \right) - \frac{\delta t}{ \alpha} g \left( (\zeta_e^l)^{m-1}_{l_0+r} -  (\zeta_e^r)^{m-1}_{l_0-r}  \right),
		\end{equation} 
		which gives $ (q_e)^m_{l_0-r}$ in terms of its values at the previous time step  and in terms of interior points. Now we can solve $(q_e)^m_{l_0-r}$ immediately. Once the value of $(q_e)^m_{l_0-r}$ is obtained, we can find the values of $\zeta^l_e |_{x=l_0-r}$ and $\zeta^r_e |_{x=l_0+r}$ by using equations \eqref{qobject} and the transport equations for the Riemann invariants as the strategy in Section \ref{discon}.

		\subsubsection{At the end of the chamber}
		The corresponding boundary condition at the end of the chamber, located at $x=l_1$, is given by
		$${q_e}_{|_{x=l_1}} = 0.$$ Hence, recalling the definition of the right-going Riemann invariant $R$, we recover the surface elevation $\zeta_e$ at $x=l_1$, namely
		$$ {\zeta_e}_{|_{x=l_1}} = \frac{1}{g}\left(\frac{R_{|_{x=l_1}}}{2} + \sqrt{g h_0}\right)^2 -h_0.$$

	\section{Numerical validations}
	In this section, we use the scheme introduced in Section \ref{discret} to simulate our model. For the fluid, we always consider the density of water $\rho = 1000 \,\mathrm{kg}/\mathrm{m}^{3}$ and the gravitational acceleration $g= 9.81\, \mathrm{m}/\mathrm{s}^{2}$. The entry of the domain is set at $x=-l= -30 \,\mathrm{m}$ and the prescribed function $f$ is given by 
	\begin{equation*}
f(t)=\sin \left(\frac{2 \pi}{T}\, t \right),
	\end{equation*}
	where $T=1.5 \, \mathrm{s}$ is the period. Using the notations as before, we consider $l_0= 11 \,\mathrm{m}$, $r=1 \,\mathrm{m}$ and $l_1= 17  \,\mathrm{m}$ and the fluid height at rest before the step $h_s=15 \, \mathrm{m}$.
	We compute the solution by using the Lax-Friedrichs scheme in the exterior domain $[-30,10]\cup [12,17]$, with a refined mesh with $N_x=2300$ and a time step $\delta_t = \frac{0.7 }{\sqrt{gh_s}}\delta_x$ with space step $\delta_x = 0.02 \, \mathrm{m}$. Here, the CFL number is $0.7$, which is commonly used to prescribe the terms of the finite-difference approximation of a PDE (see for instance \cite{ozicsik2017finite}).
In the interior domain, the solution can be computed using the transmission conditions \eqref{transwalls} with $h_w= h_0 + \zeta_w$ and $\zeta_w=\, -7.5\, \mathrm{m}$.
	
	\subsection{Numerical solutions}\label{stepwithout}
\begin{figure}
		\begin{tabular}{cc}
			\includegraphics[width=7.5cm]{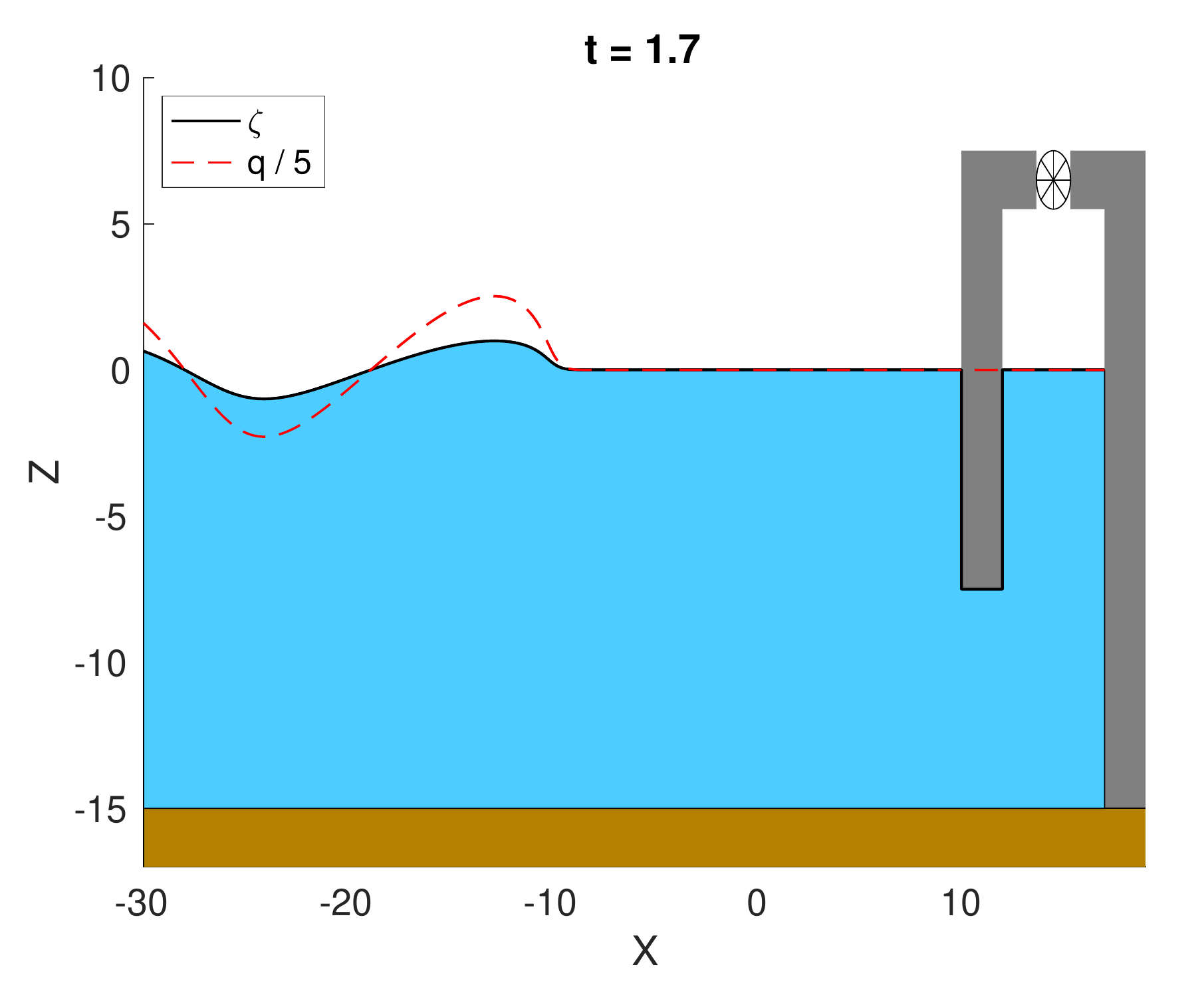}&\includegraphics[width=7.5cm]{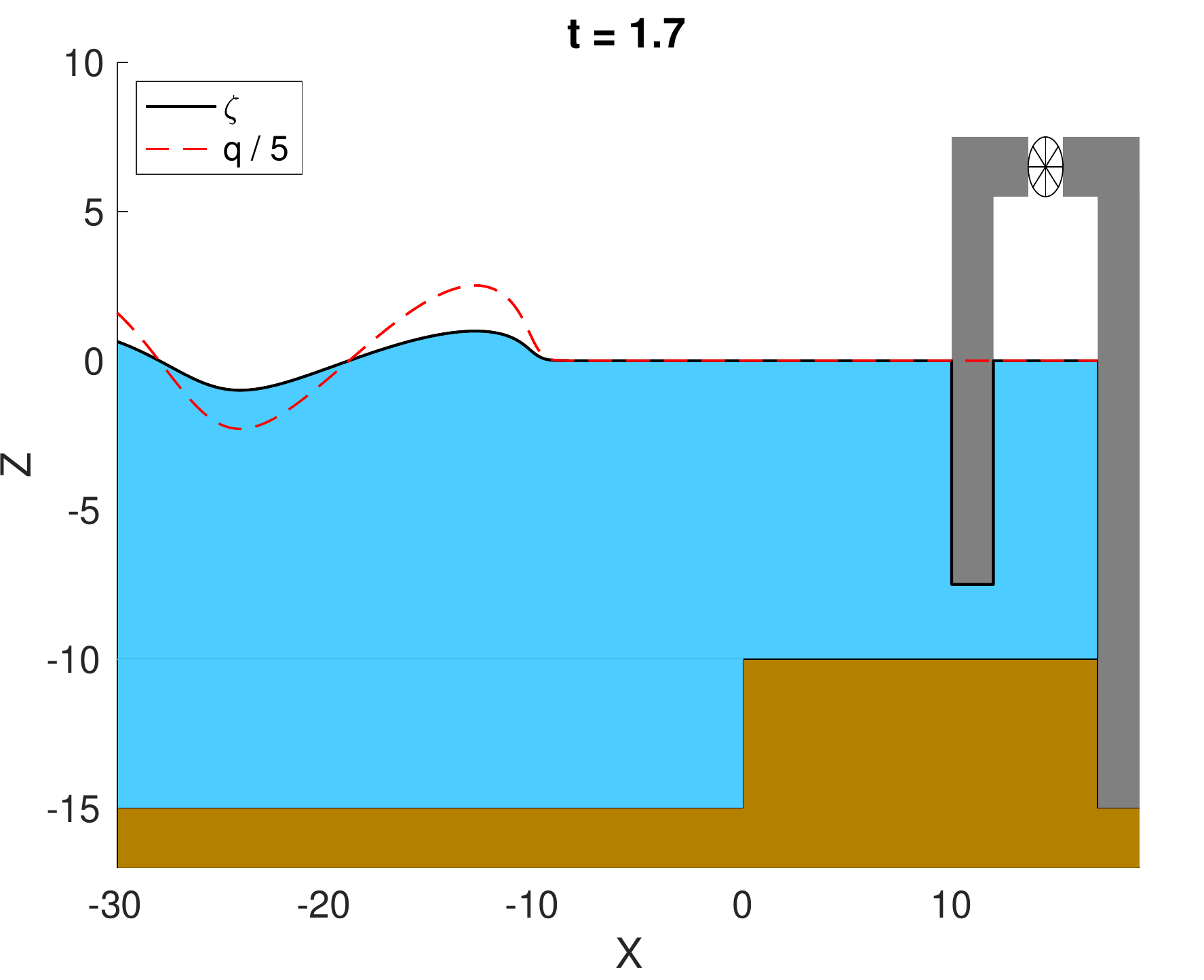}\\
			(a) & (b) \\
		\end{tabular}	

		\begin{tabular}{cc}
			\includegraphics[width=7.5cm]{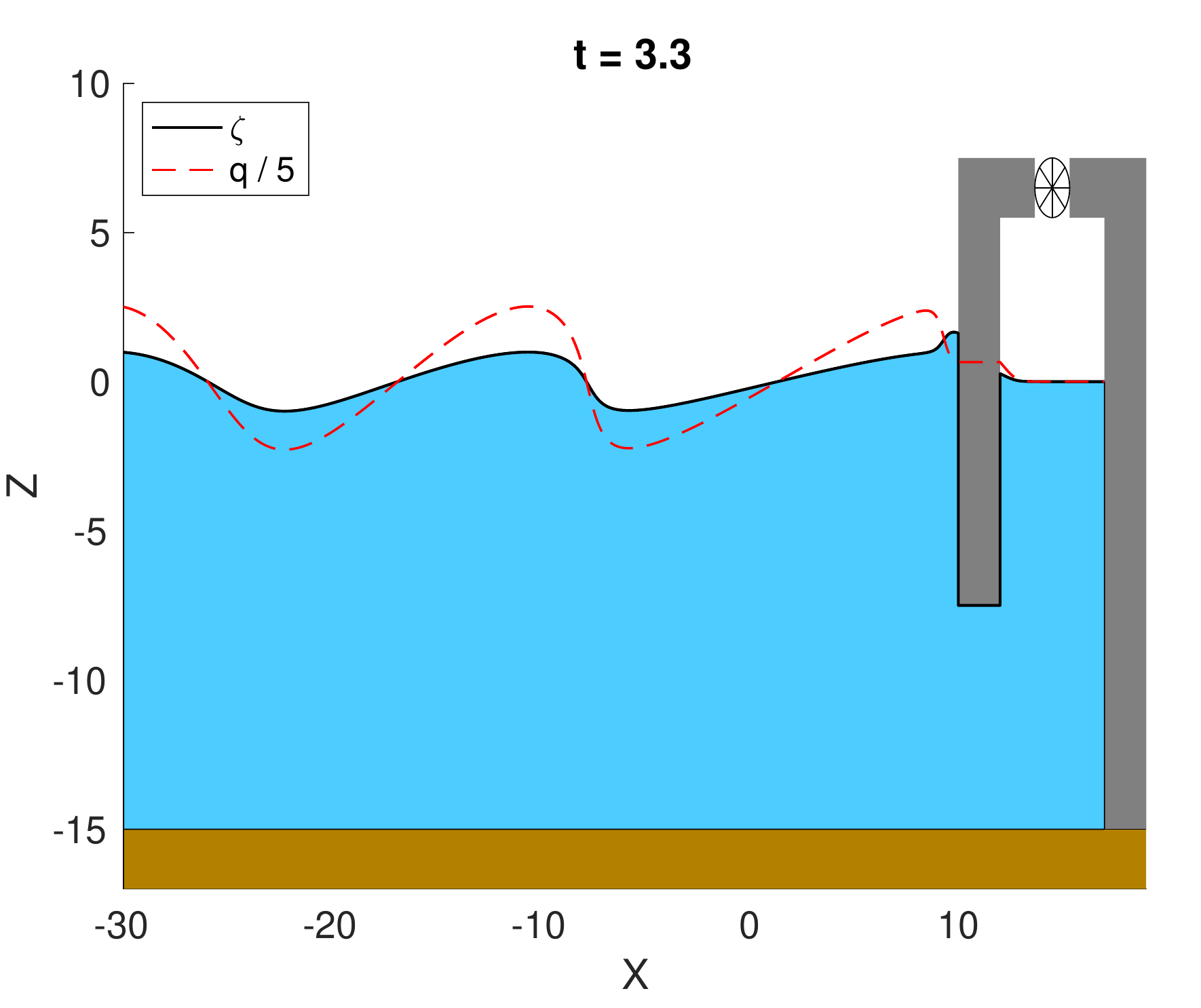}&\includegraphics[width=7.5cm]{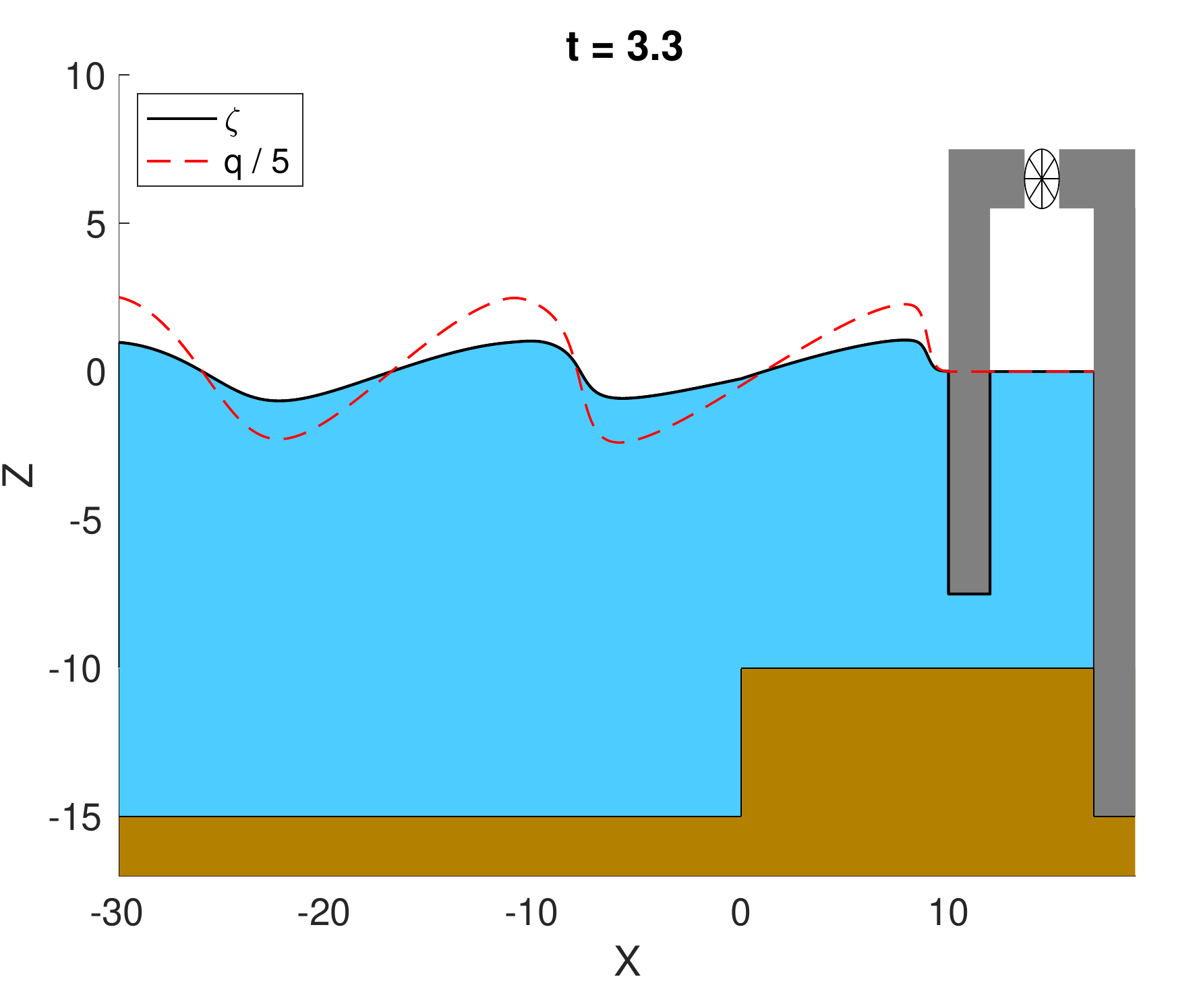}\\	
			(c) & (d) \\
		\end{tabular}	

		\begin{tabular}{cc}
			\includegraphics[width=7.5cm]{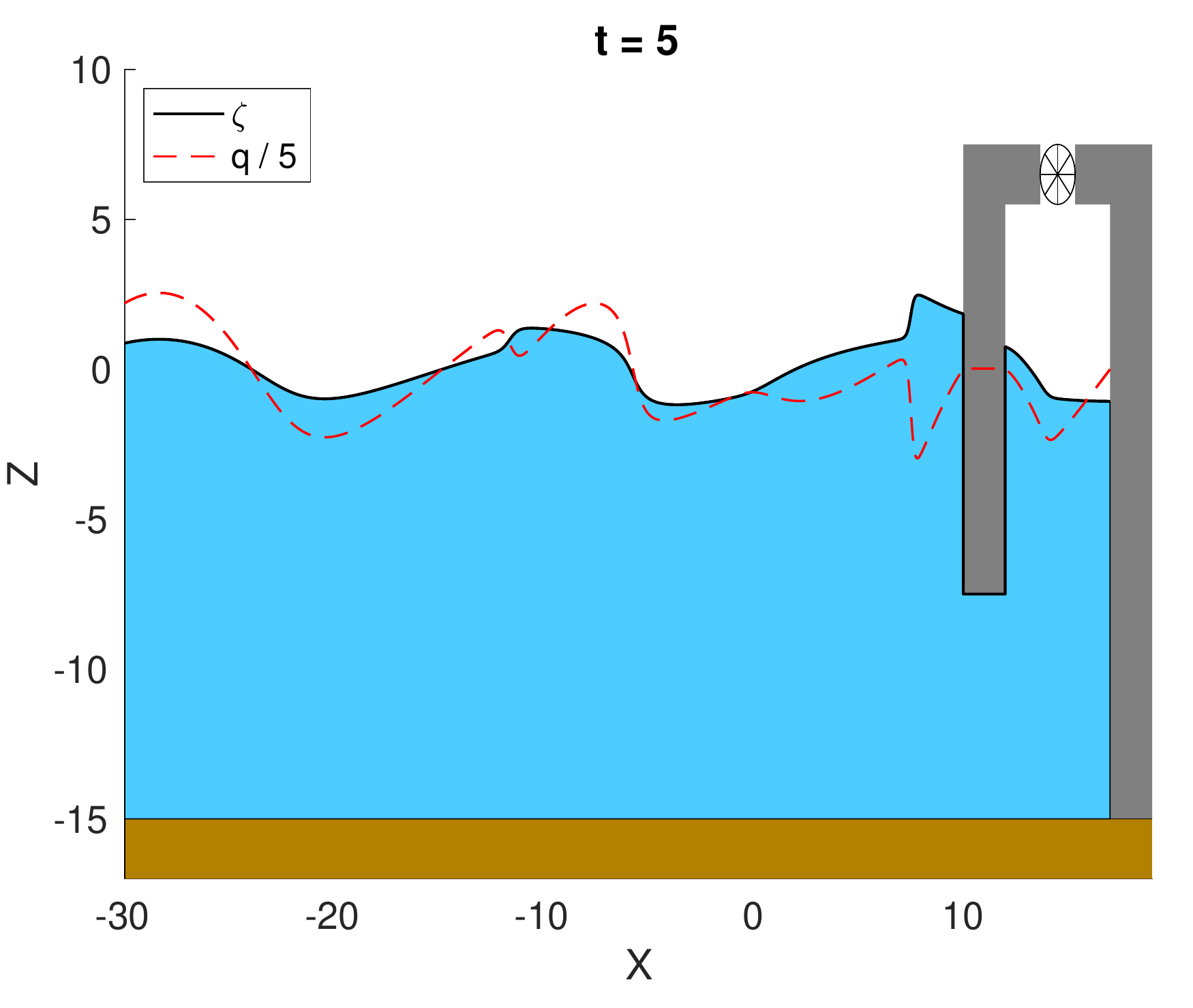}&\includegraphics[width=7.5cm]{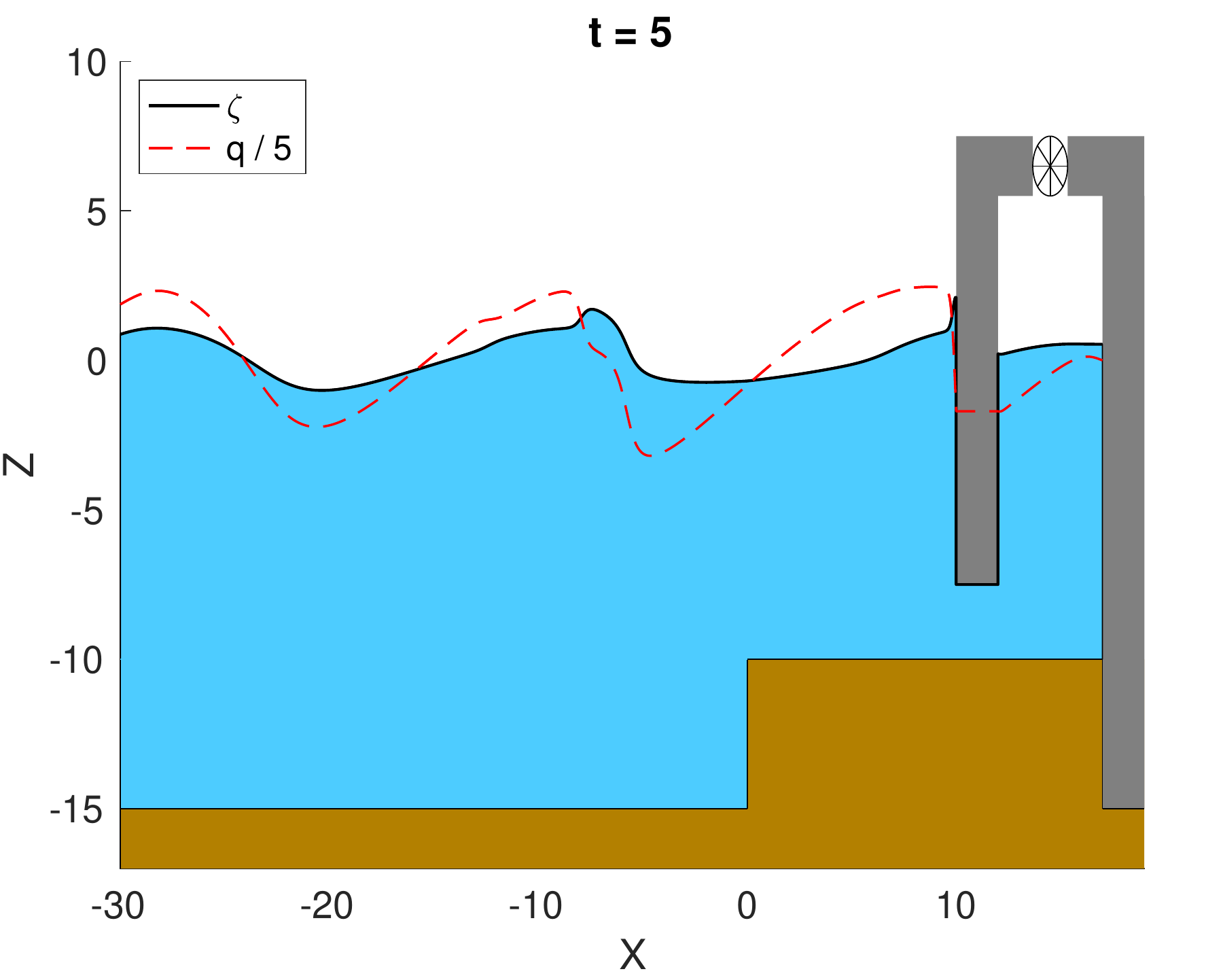}\\	
			(e) & (f)
		\end{tabular}	
		\caption{Comparisons between the numerical results without step (left) and with step (right) at times $t=1.7 \,\mathrm{s}$, $t=3.3\,\mathrm{s}$ and $t=5\,\mathrm{s}$.}
		\label{comparisonstepvsnostep}
	\end{figure}

	In real applications, an OWC device can be deployed on a stepped sea bottom in order to improve its performance. It is important then to have a good understanding of the impact of a step in the topography. Here, we test and compare the case without step $s=0 \, \mathrm{m}.$ ($h_0=15\, \mathrm{m} $) to the case with a step of height $s=5 \, \mathrm{m} $ ($h_0=10\, \mathrm{m}\ $) considering the previous physical parameters. 
	 The numerical solutions are plotted in Figure \ref{comparisonstepvsnostep} at times $t=1.7 \,\mathrm{s}$, $t=3.3 \,\mathrm{s}$ and $t=5\, \mathrm{s}$. The plots (a), (c), (e) show the solutions without stepped bottom, while the plots (b), (d), (f) show the solutions with stepped bottom. 

   We find that, before the waves encounter the step, there is no significant difference between the OWC model without stepped bottom and with stepped bottom (see (a) and (b)). But when the waves encounter the step in the bottom and arrive into the chamber, we can see that, the waves in the OWC model without stepped bottom move significantly faster than the waves in the OWC model with stepped bottom.
    In particular, at $t=3.3 \,\mathrm{s}$ the waves in the OWC model without stepped bottom has already arrived to the chamber and will begin to change the water level in the chamber, while the waves in the OWC model with stepped bottom have not reached yet and the water will rise inside the chamber later (see (c) and (d)). As the step at bottom is a sort of obstacle for the incoming wave, this phenomenon is reasonable.

	As one may expect, the incoming wave split into two parts when it touches the left wall of the partially immmersed structure. One part enters the chamber and changes the volume of the air that makes the turbine rotate. The other part is reflected and becomes an outgoing wave, as we can see in Figure \ref{comparisonstepvsnostep}. At $t=5 \,\mathrm{s}$, the reflected wave in the OWC model without stepped bottom already reaches $x=-10 \,\mathrm{m}$, while the reflected wave in the OWC model with stepped bottom has not reached $x=-10 \,\mathrm{m}$ (see (e) and (f)). This shows that the reflected waves in the OWC model with stepped bottom move slower than the waves in the OWC model without stepped bottom. 
	
	This difference can be explained by the fact that more incident wave energy is converted when a step is added. In other words, the OWC with stepped bottom would be more efficient than the one without stepped bottom, which is in agreement with the result by Rezanejad and Soares in \cite{RezSoa18}.

	\subsection{Accuracy analysis}

	In numerical validations, accuracy analysis is of importance. As we can see in Figure \ref{OWC}, the configuration of OWC device is essentially constituted from three parts:  the domain before the step in the sea bottom, the domain after the step  and the chamber. We implement our algorithms by gathering together the three parts. It is worth mentioning that one compact algorithm is also actionable.

\begin{figure}
		\begin{tabular}{cc}
			\includegraphics[width=7.5cm]{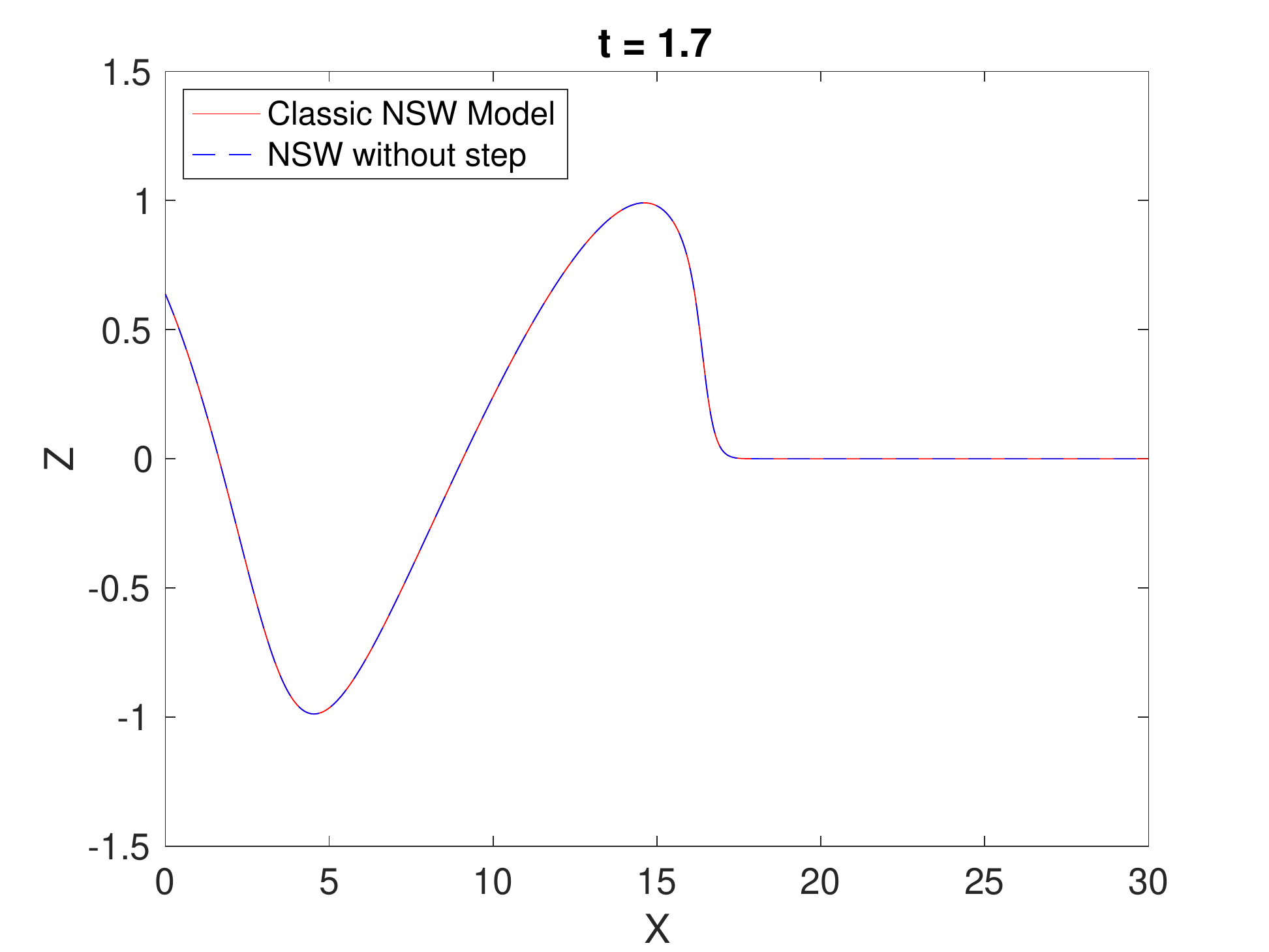}&\includegraphics[width=7.5cm]{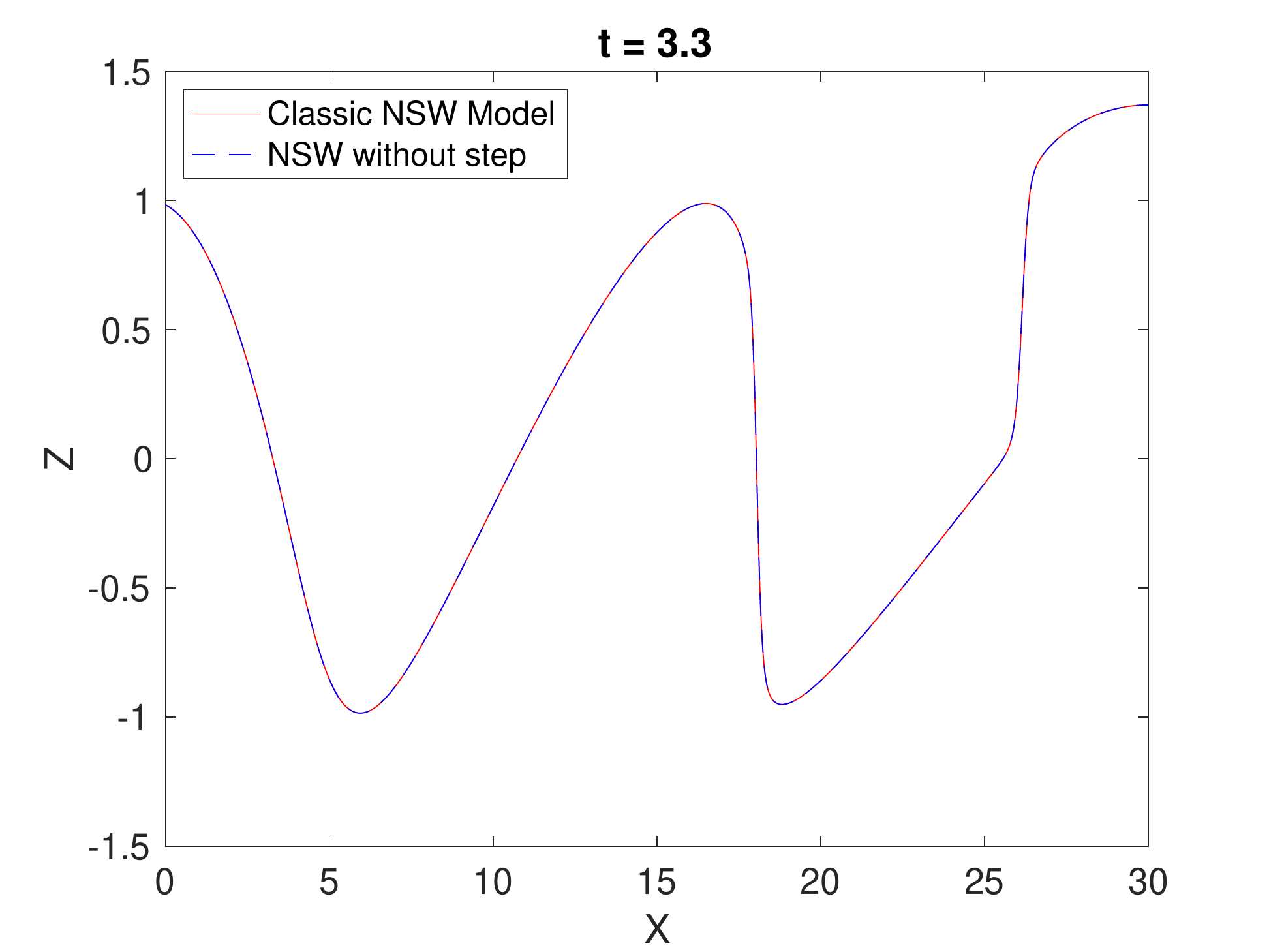}\\
			(a) & (b) 
		\end{tabular}
		\begin{tabular}{c}
			\includegraphics[width=7.5cm]{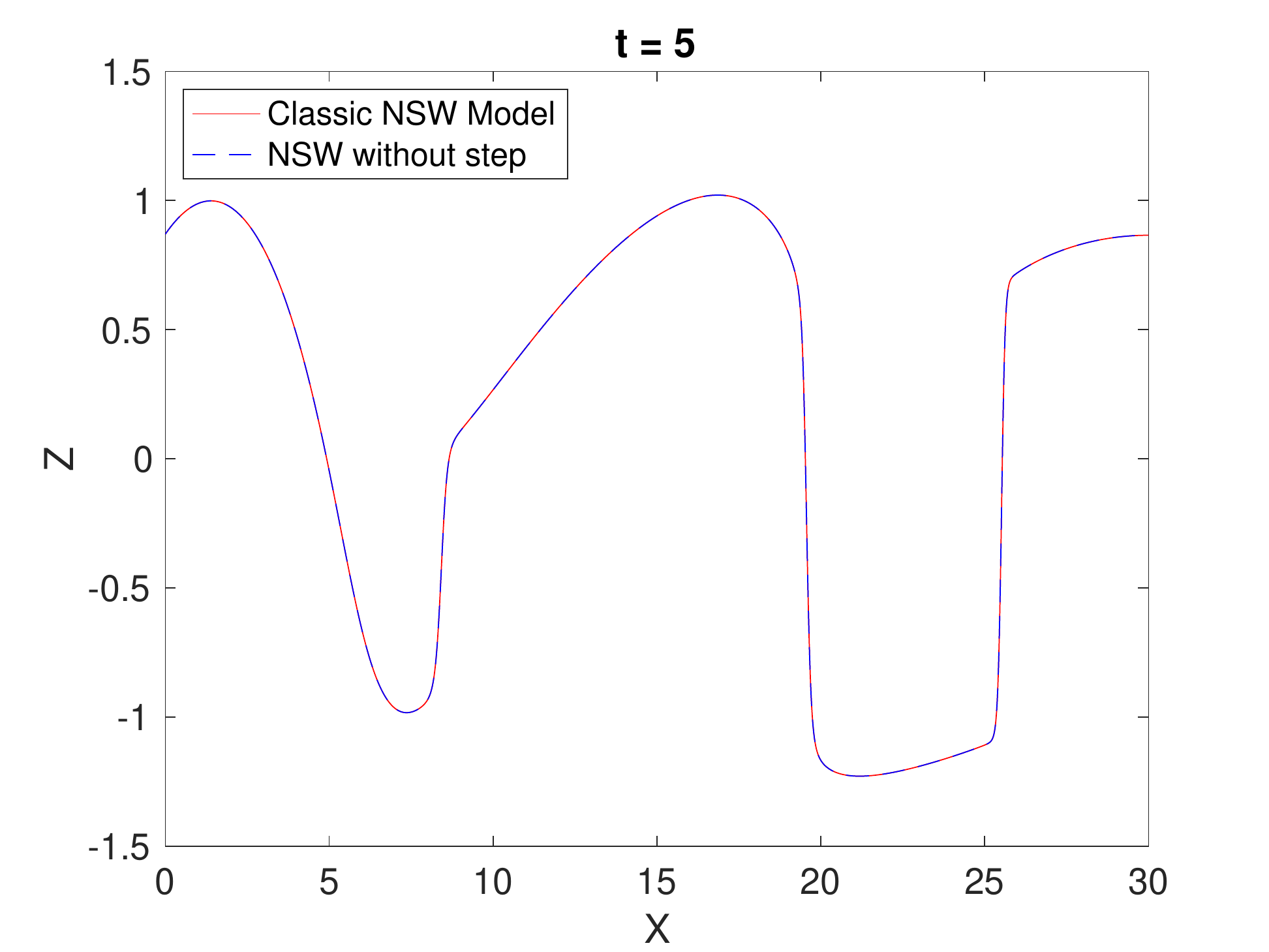}\\
			(c)  
		\end{tabular}		
		\caption{\small Comparison between classical NSW model and our model without step in different times considering
			$\delta_x = L/1500 \,\mathrm{m}$ and $L=30 \,\mathrm{m}$.}
		\label{comparison}
	\end{figure}

\begin{figure}
		\centering
		\begin{minipage}[t]{0.48\textwidth}
			\centering
			\includegraphics[width=7.5cm]{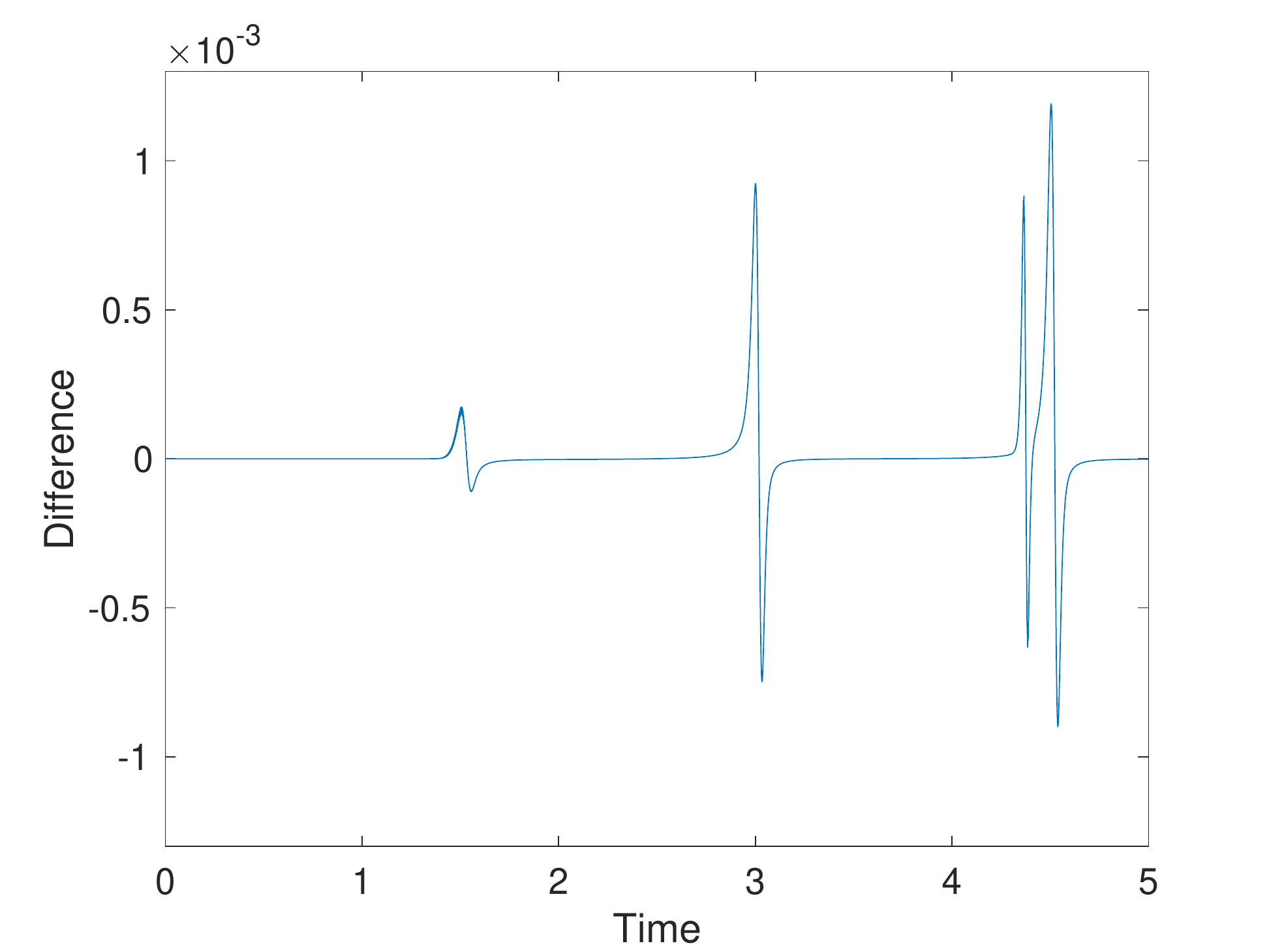}
		\end{minipage}
		\caption{\small Difference between classical NSW model and our NSW model without step considering
			$\delta_x = L/1500 \,\mathrm{m}$ and $L=30 \,\mathrm{m}$.}
		\label{comparison2}
	\end{figure}

	In order to make it possible to verify our algorithm, we do the following accuracy analysis. Under the same initial wave and physical parameters, we compare the free surface elevation $\zeta_e$ of the classical nonlinear shallow water wave model with our model without discontinuous topography. Figure \ref{comparison} shows that there is no significant difference between the two cases. 	Moreover, we also find that the error is of order $10^{-3}$ (see Figure \ref{comparison2}), which is acceptable since the Lax-Friedrichs method is first-order accurate in space.

	\subsection{Absorbed power and efficiency}\label{effpower}
	Designing a WEC of high efficiency is nowadays a hot topic in all regions and countries over the world. In this regard, we present in this section the method to calculate the absorbed power as well as the efficiency of the OWC considered in this work.
	
The primary efficiency $\eta_{Reg}$ of the device is defined by the ratio of the absorbed power from the waves to the incident wave power. 
From the seminal  work of Evans in  \cite{evans1982wave}, we know that in the linear time-harmonic theory the volume flux $Q(t)=\mathrm{Re}\{qe^{-i\omega t}\}$ is assumed linearly proportional to the pressure in the chamber $P(t)= \mathrm{Re}\{pe^{-i\omega t}\}$. Using this assumption, the average power absorbed from regular waves over one wave period, denoted by 
$P_{Reg}$, is given by
	\begin{equation}P_{Reg} = \frac{1}{2} \lambda |p|^2,
	\label{harmonicPower}
	\end{equation}
	where $p$ is the time independent and $\lambda$ is a positive constant associated with linear air turbine characteristics. On the other hand, following \cite{RezSoa18} in experiments the average power absorbed from regular waves can be determined by:
	\begin{equation}
		P_{Reg}=\frac{1}{T}\int_0^T PQ dt,
		\label{testPower}
	\end{equation}
	where $T$ is the duration of the test.
	The incident wave power $P_{inc}$ is defined as the product of total energy per wave period $E_{inc}$ and the group velocity $c_g$ (see \cite{dean1991water}):  
	\begin{equation*}
	{P_{inc}=E_{inc} \,c_{g}}, 
	\end{equation*}
	with 
	$$ {E_{inc}=\frac{1}{2} \rho g L A^{2}}, \qquad  {c_{g}= \frac{\omega}{2k}\left(1+\frac{2kh_s}{\sinh(2kh_s)}\right)},$$ 
	 and the dispersion relation given by $$\omega^2=g k \tanh(kh_s),\vspace{0.5em}$$
	 where $\omega$ is the frequency, $k$ is the wave number, $h_s$ is the fluid height at rest before the step, $\rho$ is the density of the fluid, $g$ is the gravitational acceleration and $L$ is the projected width of the WEC perpendicular to the incident
	wave direction, $A$ is the amplitude of the wave. In the shallow water regime 
	$kh_s\ll1$ and the group velocity reduces to $c_g= \sqrt{gh_s}.$
	 Thus, the primary efficiency of the device in regular wave is given by 
		$$\eta_{Reg} = \frac{P_{Reg}}{P_{inc}}.$$
		
We notice that in both \eqref{harmonicPower} and \eqref{testPower} the absorbed power (hence the primary efficiency) strongly depends on the air pressure in the chamber. In our model, it is considered to be a constant, namely the atmospheric pressure $P_{\mathrm{atm}}$. However, when the waves arrive into the chamber and change the volume of the air, the air pressure in the chamber will certainly change as well. In this case, the pressure will no more be a constant, but depends on time. Hence, to study more rigorously the absorbed power and the primary efficiency of the OWC, this fact must be taken into account in the model. This will be addressed in our future work. Analogously, the improvement of the efficiency of an OWC device deployed on a stepped sea bottom can be also investigated with a better knowledge of the air pressure in the chamber. From the results in Section \ref{stepwithout}, we can expect that significant improvements in the efficiency can be achieved by adding a step at the bottom of the sea. 
	
	\section*{Acknowledgements}
	The authors warmly thank David Lannes for his helpful comments and advises and also the organisers of the summer school CEMRACS 2019 during which this work was done. E.B. is supported by the Starting Grant project “Analysis of moving incompressible fluid interfaces” (H2020-EU.1.1.-639227) operated by the European Research Council.
	J. H. is supported by the PostDoc program Sophie
Germain of the Fondation Mathématique Jacques Hadamard.
	G. V-H. is supported by the European Union's Horizon 2020 research and innovation programme under the Marie Sklodowska-Curie grant agreement No 765579.

	
	\nocite{*} 
	\bibliographystyle{siam}
	\bibliography{biblio}
	
\end{document}